\documentclass[fleqn,12pt]{article}
\usepackage[english]{babel}
\usepackage[latin1]{inputenc}
\usepackage{graphicx,amsmath,amsfonts,amssymb,amscd,bezier,amstext,makeidx}
\usepackage{url}

\newtheorem{theorem}{Special Theorem}[section]

\newtheorem{theo}[theorem]{Theorem}
\newtheorem{lemm}[theorem]{Lemma}
\newtheorem{corol}[theorem]{Corollary}
\newtheorem{remark}[theorem]{Remark}

\newtheorem{prop}[theorem]{Proposition}
\newtheorem{example}[theorem]{Example}

\title{On a class of Rauzy fractals without the finiteness property}
\author{Gustavo A. Pavani }
\date{November, 2016}

\begin{document}

\maketitle

\begin{abstract}
\noindent We present some topological and arithmetical aspects of a class of Rauzy fractals $\mathcal{R}_{a,b}$ related to the polynomials of the form $P_{a,b}(x)=x^{3}-ax^{2}-bx-1$, where $a$ and $b$ are integers satisfying $-a+1 \leq b \leq -2$. This class has the property that $0$ lies on the boundary of $\mathcal{R}_{a,b}$. We construct  explicit finite automata that recognize the boundaries of these fractals, which allows to establish the number of neighbors of $\mathcal{R}_{a,b}$. In particular,  we prove that if $2a+3b+4 \leq 0$ then $\mathcal{R}_{a,b}$ is not homeomorphic to a topological disk. We also show that the boundary of the set $\mathcal{R}_{3,-2}$ is generated by two infinite iterated function systems.
\end{abstract}

\section{Introduction}

Our aim is to study a class of Rauzy fractals related to an algebraic integer $\beta$ which does not satisfy a certain property called (F) Property. This study envolves fractal tilings, automata, $\beta-$numeration, and infinite iterated function systems (IIFS).\\

The Rauzy fractal was introduced by G. Rauzy in 1982 \cite{Rauzy1982} and it is the set

\centerline{$\displaystyle \mathcal{E}=\left\{\sum_{i=0}^{+\infty}\ell_{i}\alpha^{i},\, \ell_{i}\in\{0,1\},\, \ell_{i}\ell_{i+1}\ell_{i+2}=0, \, \forall i \geq 0 \right \}$,}

where $\alpha$ is one of the complex roots of the polynomial $x^{3}-x^{2}-x-1$. 

There are several ways to construct the Rauzy fractal, one of them is by using substitutions. A substitution $\sigma$ is a map from an alphabet $\cal{A}$ into the set $\cal{A}^{*}$ of nonempty finite words on $\cal{A}^{*}$. The initial aim of Rauzy was to establish a geometric representation to the symbolic dynamical system associated with the substitution $\sigma$ given by $\sigma(0)\rightarrow 01$, $\sigma(1)\rightarrow 02$ $\sigma(2)\rightarrow 0$. Since then, this set and its generalizations have been studied by many mathematicians, due to its strong connections with other fields of mathematics, for instance, tilings \cite{Praggastis, Akiyama2000, ArnouxIto}, numeration systems \cite{Me2000, Me2006, Praggastis}, Markov partitions for toral automorphisms \cite{Me98, Praggastis, Vershik}, geometric representation of symbolic dynamical systems \cite{CanteriniSiegel, ArnouxIto, Sano, ArnouxRauzy, Me98, Sirvent, Zamboni, SiegelThus}, simultaneous diophantine approximations \cite{Arnoux2002, Chekhova, HubMess}, and the theory of quasicrystals \cite{QuasiCrystals}.

The Rauzy fractal has remarkable properties: it is a compact and connected subset of $\mathbb{C}$, its interior is simply connected, and it induces a periodic tiling of the complex plane modulo $\mathbb{Z}\alpha^{-3}+\mathbb{Z}\alpha^{-2}$. Moreover, it is divided into three self-similar copies of itself which correspond to an exchange of domains \cite{Rauzy1982}. 

Another way to obtain the Rauzy fractal is via $\beta-$representation. Given a real number $\beta>1$, a $\beta$-representation (or $\beta$-expansion) of a number $x \in \mathbb{R}^{+}$ is an infinite sequence $(x_{i})_{i\leq k}$, where $k \in \mathbb{Z}$, $x_{i} \geq 0$ such that $x=\sum_{i=-\infty}^{k}x_{i}\beta^{i}$. The digits $x_{i}$  can be computed using the greedy algorithm as follows (see \cite{Parry, Frougny2000} for details): denote by $\left\lfloor x\right\rfloor$ and $\left\{x\right\}$ the integer and fractional parts of the number $x$. There exists $k \in \mathbb{Z}$ such that $\beta^{k} \leq x < \beta^{k+1}$. Let $x_{k}=\left\lfloor x/\beta^{k}\right\rfloor$ and $q_{k}=\left\{x/\beta^{k}\right\}$. Then, for $i<k$, put $x_{i}=\left\lfloor \beta q_{i+1}\right\rfloor$ and $q_{i}=\left\{\beta q_{i+1}\right\}$. We obtain $x=x_{k}\beta^{k}+x_{k-1}\beta^{k-1}+\cdots$. If $k<0$ ($x<1$) we put $x_{0}=x_{1}=\cdots=x_{k+1}=0$. If a $\beta$-representation ends with infinitely many zeros, it is said to be finite and the ending zeros can be omitted. Then, the sequence will be denoted by $(x_{i})_{n\leq i \leq k}$ or $x_{k}\cdots x_{n}$. If $\beta$ is an integer, the digits $x_{i}$ belong to the set $B=\{0, \cdots, \beta\}$, or to the set $B=\{0, \cdots, \left\lfloor \beta\right\rfloor\}$, otherwise. 

In particular, when $\beta$ is a Pisot number, i.e., an algebraic integer greater then 1 such that all its Galois conjugates have modulus less than 1, we obtain classes of Rauzy fractals associated to these Pisot numbers. Cubic and unitary Pisot numbers were classified by Akiyama in \cite{Akiyama2000} as being exactly the set of dominant roots of the polynomial $P_{a,b}(x)=x^3-ax^{2}-bx-1$, satisfying one of the following conditions

\hspace{0.50cm} \textbf{a)} $1 \leq b \leq a$ and $d(1,\beta)=ab1$;

\hspace{0.50cm} \textbf{b)} $b=-1$, $a \geq 2$ and $d(1,\beta)=.(a-1)(a-1)01$;

\hspace{0.50cm} \textbf{c)} $b=a+1$ and $d(1,\beta)=.(a+1)00a1$;

\hspace{0.50cm} \textbf{d)} $-a+1 \leq b \leq -2$ and $d(1,\beta)=.(a-1)(a+b-1)(a+b)^{\infty}$,

where $(a+b)^{\infty}$ is the periodic expansion $(a+b)(a+b)(a+b)\ldots,$ and $d(1,\beta)$ is the R\'{e}nyi $\beta$-representation of 1 (see \cite{Renyi} for the definition). \\

Let Fin($\beta$) be the set of nonnegative real numbers that have a finite $\beta$-representation. We say that a Pisot number $\beta$ has the finiteness property (or (F) property ) if $\mathbb{Z}[\beta] \cap [0,+\infty[ \subset$ Fin($\beta$). Therefore, the Pisot numbers in the sets \textbf{a)}, \textbf{b)} and \textbf{c)} has the (F) property, while the Pisot numbers in \textbf{d)} have not. Many works were done for the classes \textbf{a)} and \textbf{b)} (see \cite{Thus2006, ItoKimura, Me2000, Me2013, Me2015, Loridant}). \\

In this paper we will study the properties of the Rauzy fractals associated to the class of Pisot numbers which do not satisfy the (F) Property, that is, the case where $-a+1 \leq b \leq -2 $. As we shall see, this class shares common features with the others previously studied. For instance, these fractals sets are compact and they tile the plane. In fact, this class of fractals can be obtained via $\beta-$substitution defined by $\sigma(1) \mapsto 1^{(a-1)}2$, $\sigma(2) \mapsto 1^{(a+b-1)}3$, $\sigma(3) \mapsto 1^{(a+b)}3,$ (see \cite{BertheSiegel}) and some topological properties for fractal sets arising from Pisot substitutions are known (see \cite{Canterini}). On the other hand, zero is not an inner point for the fractals of this class, as it occurs in the classes for which the (F) property holds. In this work, we obtain explicit finite state automata that generate the boundary of $\mathcal{R}_{a,b}$. These automata lead to several results: we obtain a formula for the number of the neighbors of $\mathcal{R}_{a,b}$ in the periodic tiling and we prove that if $2a+3b+4 \leq 0$, then $\mathcal{R}_{a,b}$ is not homeomorphic to a topological disk. We study in more details the boundary of the set $\mathcal{R}_{3,-2}$, in particular we prove that the boundary of $\mathcal{R}_{3, -2}$ is generated by two infinite iterated function systems. Notice that the set $\mathcal{R}_{3,-2}$ is related with Special Pisot numbers, i.e., a Pisot number $\beta$ such that $\beta/(\beta-1)$ is also a Pisot number (see \cite{LagariasI, Smith}). \\
 
This paper is organized in this way. In the second section we briefly describe the $\beta-$numeration necessary to define the Rauzy fractal sets that we are considering. In the third section we give some properties of the boundary of $\mathcal{R}_{a,b}$. In the fourth section we construct the automata that recognize the boundaries of the sets $\mathcal{R}_{a,b}$, and in fifth section we use an automaton to obtain two IIFS for the boundary of $\mathcal{R}_{3,-2}$, and we show a geometric method that could be used for parametrizing the boundary of $\mathcal{R}_{3,-2}$.  

\section{Numeration system and Rauzy fractal}

In the sequel, we will suppose that $\beta$ is a cubic and unitary Pisot number which does not satisfy the (F) property and we will denote by $\alpha$ and $\lambda$ its Galois conjugates. Let $P_{a,b}(x)=x^{3}-ax^{2}-bx-1$ be the minimal polynomial of $\beta$. Next, we will consider a generalization of numeration system induced by the $\beta$-expansion which only can be applied on integer numbers.\\

Let $(T_{n})_{n\geq0}$ be the recurrent sequence defined by $T_{0}=1$, $T_{1}=a$, $T_{2}=a^{2}+b, T_{n+3}=aT_{n+2}+bT_{n+1}+T_{n}$, satisfying the condition $-a+1 \leq b \leq -2$ for all $ n \geq 0$.

\begin{prop} \label{todoint} Every nonnegative integer $n$ can be uniquely expressed as $n=\sum_{i=0}^{N}\ell_{i}T_{i}$, where $\ell_{i} \in \{0, \ldots, a-1\}$ and $\ell_{j}\ell_{j-1}\cdots \ell_{j-k} \leq_{lex}(a-1)(a+b-1)(a+b)\cdots (a+b)$, for all $j \geq k \geq 0$, where ``$\leq_{lex}$'' is the lexicographical order. \end{prop}

For the proof we need the following Lemma:

\begin{lemm} The sequence $(T_{n})_{n \geq 4}$ satisfies

\begin{center} $T_{n}=(a-1)T_{n-1}+ (a+b-1)T_{n-2}+(a+b)T_{n-3}+\cdots +(a+b)T_{1}+(a+b+1)T_{0}$,  \end{center}

for all $n \geq 4$. \end{lemm} 

\noindent \textbf{Proof.} The proof is by recurrence on $n$. It is not difficult to verify that the relation is valid for $n=4,5,6$. Suppose that the relation holds for all $k<n$. Since $T_{n}=aT_{n-1}+bT_{n-2}+T_{n-3}$, then $ T_{n}=aT_{n-1}+bT_{n-2}+Q$, where

$Q=(a-1)T_{n-4}+(a+b-1)T_{n-5}+ (a+b)T_{n-6}+\cdots+(a+b)T_{1}+(a+b+1)T_{0}$. 

Then,

$\begin{array}{ll}
  aT_{n-1}+bT_{n-2}+(a-1)T_{n-4}+(a+b-1)T_{n-5}  & = (a-1)T_{n-1}+(a+b-1)T_{n-2} \\
   & +(a+b)T_{n-3}+(a+b)T_{n-4} \\
	& + (a+b)T_{n-5} \\
   
\end{array}$ \\

In fact,

$aT_{n-1}+bT_{n-2}+(a-1)T_{n-4}+(a+b-1)T_{n-5}=$

\hspace{1.15cm}$=(a-1)T_{n-1}+(a+b)T_{n-2}+bT_{n-3}+aT_{n-4}+(a+b-1)T_{n-5}$

\hspace{1.15cm}$ =(a-1)T_{n-1}+(a+b-1)T_{n-2}+(a+b)T_{n-3}+(a+b)T_{n-4}+(a+b)T_{n-5}.$  $\Box$\\

\textbf{Proof of Proposition \ref{todoint}.} The digits $(\ell_{j})_{0 \leq j \leq k(N)}$ are obtained by using the greedy algorithm. Since  $-a+1 \leq b \leq -2$, then $(T_{n})_{n \geq 0}$ is an increasing sequence of natural integers. Hence, by the  definition of the greedy algorithm, we can prove that \\

\centerline{$\sum_{i=0}^{j}\ell_{i}T_{i}<T_{j+1}$,}

\bigskip

for all $0 \leq j \leq k(N)$  (see \cite{Parry}). Thus, 

\begin{equation} \label{digits} \ell_{j}\ell_{j-1} \cdots \ell_{j-k} <_{lex} (a-1)(a+b-1)(a+b) \cdots (a+b)(a+b+1), \forall j \geq k \geq 0. 
\end{equation}

Therefore, by using (\ref{digits}), we obtain that $\ell_{j} \cdots \ell_{j-k}\leq_{lex}(a-1)(a+b-1)(a+b) \cdots (a+b)$. $\Box$\\

Let $\mathcal{L}=\{(\ell_{i})_{i\geq k}, k\in \mathbb{Z}, \forall n\geq k, \ell_{n}\cdots \ell_{n-k}\leq_{lex}(a-1)(a+b-1)(a+b)\cdots (a+b)\}$. Then, the Rauzy fractal is the set \\

\centerline{$\displaystyle \mathcal{R}:=\mathcal{R}_{a,b}=\{\sum_{i=2}^{+\infty}\ell_{i}\theta_{i}, \, (\ell_{n})_{n\in \mathbb{Z}} \in \mathcal{L}\}$}

\bigskip

where $\theta_{i}=\alpha^{i}$, if $\alpha \in \mathbb{C}\setminus \mathbb{R}$ or $\theta_{i}=(\alpha^{i},\lambda^{i})$, if $\alpha \in \mathbb{R}$. Observe that  $\mathcal{R} \subset \mathbb{C}$ or $\mathcal{R} \subset \mathbb{R}^{2}$.

\begin{remark} We take the summation beginning from 2 in the definition of the Rauzy fractal just for technical purposes. 
\end{remark}

\begin{example}

1. If $a=3$ and $b=-2$, we can show that $P_{3,-2}(x)=x^{3}-3x^{2}+2x-1$ has one real root $\beta>1$ and two complex conjugates roots $\alpha, \overline{\alpha}$ wich satisfy $|\alpha|,|\overline{\alpha}|<1$. In this case $\alpha \approx 0.33764+0.56228i$. The Rauzy fractal (Figure 1) is

\begin{center}$\mathcal{R}_{3,-2}=\{\sum_{i=2}^{+\infty}\ell_{i}\alpha^{i}, \forall j\geq n \geq 2, \ell_{j}\ell_{j-1}\cdots \ell_{n} \leq_{lex}201\cdots 1\}.$ \end{center}

2. If $a=6$ and $b=-5$, then $P_{6,-5}(x)=x^{3}-6x^{2}+5x-1$ has three real roots: $\beta \approx 5.048917340$, $\alpha \approx 0.3079785280$ and $\lambda \approx 0.6431041320$. The Rauzy fractal (Figure 2) in this case is 

\begin{center}$\mathcal{R}_{6,-5}=\{(\sum_{i=2}^{+\infty}\ell_{i}\alpha^{i},\sum_{i=2}^{+\infty}\ell_{i}\lambda^{i}), \forall j\geq n \geq 2, \ell_{j}\ell_{j-1}\cdots \ell_{n} \leq_{lex}501\cdots 1\}.$ \end{center}

\end{example}

\begin{figure}[h!]
\centering
 \includegraphics[scale=0.4]{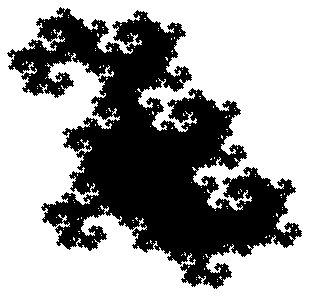}
 \includegraphics[scale=.4]{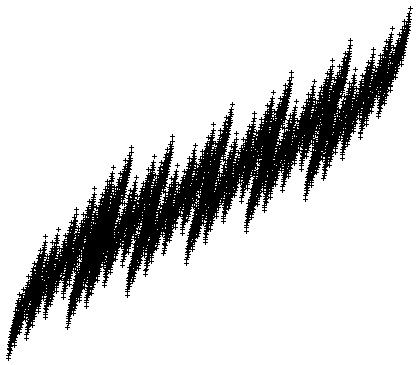}
\caption{The set $\mathcal{R}_{3,-2}$. \hspace{2.0cm} Figure 2: The set $\mathcal{R}_{6,-5}$}
\end{figure}

\section{Boundary of $\mathcal{R}$}

In this section we show some properties concerning the boundary of the Rauzy fractals. We will denote the interior of the set $\cal{R}$ by $int(\cal{R})$. We have the following Theorem (see Figures 3 and 4):

\begin{theo} \label{tiling} The Rauzy fractal induces a periodic tiling of the plane $\mathcal{P}$ modulo group $G$, that is, $\mathcal{P}= \bigcup_{u \in G}$($\cal{R}$$+u$) and int($\cal{R}$$+u$) $\cap$ ($\cal{R}$$+v$)$\neq \emptyset$ implies that $u = v$. When $\alpha \in \mathbb{C}\setminus \mathbb{R}$, then $\mathcal{P}=\mathbb{C}$ and $G=\mathbb{Z}\alpha^{-2}+\mathbb{Z}\alpha^{-1}$. If $\alpha \in \mathbb{R}$, then $\mathcal{P}=\mathbb{R}^{2}$ and $G=\mathbb{Z}(1,1)+\mathbb{Z}(\alpha,\theta)$.
\end{theo}

\noindent \textbf{Proof.} The proof can be deduced from the work of Rauzy, (see also \cite{Me2015, CanteriniSiegel}), which is done for the case when the associated $\beta-$expansion is finite, but it runs in the same way to the case we are treating. $\Box$

\begin{prop} \label{frontx} The boundary $\partial \mathcal{R}$ of $\mathcal{R}$ satisfies the property:

$\partial \mathcal{R}$ $=\bigcup_{u \in H}$$\mathcal{R}$ $\cap$ ($\mathcal{R}$$+u$), where $H$ is finite subset of $\mathbb{Z}\alpha^{-2}+\mathbb{Z}\alpha^{-1}$, whose cardinality is even and greater then or equal to 6, and $\{\pm (1+(b+1)\alpha), \pm \alpha,\pm (1+b\alpha)\} \subset H$.
\end{prop}

For the proof we need the following result (see Figure 7): 

\begin{lemm} \label{pontow} Let $\psi:\{0,1,\cdots,a-1\}^{\mathbb{N}}\rightarrow \mathbb{C}$ defined by $\psi(\ell_{0}\ell_{1}\cdots)=\sum_{i=0}^{\infty}\ell_{i}\alpha^{i}$. Let 

\hspace{1.0cm} $w_{1}  = \psi(0000((b+2)(a+b)(a-2))^{\infty})$, 

\hspace{1.0cm} $w_{2}  = \psi(01b(a-1)000((b+2)(a+b)(a-2))^{\infty})$, 

\hspace{1.0cm} $w_{3}  =  \psi(1(b+1)(a+b)(a-2)000((b+2)(a+b)(a-2))^{\infty})$, 

\hspace{1.0cm} $z_{1}=\psi(1b(a-1))$, $z_{2}=\psi(000(b+2)(a+b+1)(a+b)^{\infty}).$ 

Then $w_{1}=w_{2}=w_{3}$, $z_{1}=z_{2}$, and hence  $w_{1} \in \mathcal{R} \cap (\mathcal{R}+\alpha) \cap (\mathcal{R}+1+(b+1)\alpha)$, and also $z_{1} \in \mathcal{R} \cap (\mathcal{R}+1+b\alpha)$.
\end{lemm}

\textbf{Proof}. Let us show that $w_{1}=w_{2}$. We have,\\

 \centerline{$w_{1}  = \frac{1}{1-\alpha^{6}}((b+2)\alpha^{4}+(a+b)\alpha^{5}+(a-2)\alpha^{6})$}

 and 

\centerline{$w_{2}  =  \alpha+b\alpha^{2}+(a-1)\alpha^{3}+ \frac{1}{1-\alpha^{6}}((b+2)\alpha^{7}+(a+b)\alpha^{8}+(a-2)\alpha^{9})$. }

Then,

$\begin{array}{ll}
  w_{1}-w_{2}=0 & \Longleftrightarrow (b+2)\frac{(\alpha^{4}-\alpha^{7})}{1-\alpha^{6}}+(a+b)\frac{(\alpha^{5}-\alpha^{8})}{1-\alpha^{6}}+(a-2)\frac{(\alpha^{6}-\alpha^{9})}{1-\alpha^{6}}\\
	              & \hspace{.85cm}-\alpha-b\alpha^{2}-(a-1)\alpha^{3}=0 \\
   & \Longleftrightarrow \alpha^{4}(b+2)\frac{(1-\alpha^{3})}{1-\alpha^{6}}+\alpha^{5}(a+b)\frac{(1-\alpha^{3})}{1-\alpha^{6}}+\alpha^{6}(a-2)\frac{(1-\alpha^{3})}{1-\alpha^{6}}\\
	& \hspace{.85cm}-\alpha-b\alpha^{2}-(a-1)\alpha^{3}=0 \\
   
\end{array}$ \\

Multiplying the last equation by $(1+\alpha^{3})$ we obtain

\centerline{$\alpha^{4}(b+2)+\alpha^{5}(a+b)+\alpha^{6}(a-2)-\alpha(1+\alpha^{3})-b\alpha^{2}(1+\alpha^{3})-\alpha^{3}(a-1)(1+\alpha^{3})=0.$}

Now, developping the left side of the above equation and using the fact that $\alpha^{3}=a\alpha^{2}+b\alpha+1$ we obtain that $w_{1}=w_{2}$. 

The other cases, left to the reader, can be done in the same way. $\Box$\\

\noindent \textbf{Proof of Proposition \ref{frontx}.} Let $z \in \partial \mathcal{R}$. Since $\mathbb{C}= \bigcup_{u \in \mathbb{Z}\alpha^{-2}+\mathbb{Z}\alpha^{-1}} (\mathcal{R}+u)$ and $\mathcal{R}$ is closed, then there exists a sequence $(z_{n})_{n \geq 0}$ of elements of $\mathbb{C}$ such that

\bigskip

\centerline{lim $z_{n}=z$ and $z_{n} \notin$ $\mathcal{R}$, $\forall n\geq 0$.}

\bigskip

Then, by Theorem \ref{tiling} there exists a sequence $(u_{n})_{n \geq 0}$ of elements of $\mathbb{Z}\alpha^{-2}+\mathbb{Z}\alpha^{-1}\backslash \{0\}$ such that $z_{n} \in$ $\mathcal{R}+u_{n}$, for all $ n\geq 0$. Hence $(u_{n})_{n \geq 0}$ is bounded. Since $\mathbb{Z}\alpha^{-2}+\mathbb{Z}\alpha^{-1}$ is a lattice, then $(u_{n})_{n \geq 0}$ is a sequence that have a finite number of terms. Thus, there exists a sub-sequence $(u_{k_{n}})_{n \geq 0}$ of $(u_{n})_{n \geq 0}$ such that $u_{k_{n}}=u \in \mathbb{Z}\alpha^{-2}+\mathbb{Z}\alpha^{-1} \backslash \{0\}$. Since $z_{k_{n}} \in \mathcal{R}+u_{k_{n}}=\mathcal{R}+u$ we have $z=$lim $z_{k_{n}} \in \mathcal{R}+u$, because $\mathcal{R}+u$ is a closed set. Hence, $\partial \mathcal{R} \subset \bigcup_{u \in \mathbb{Z}\alpha^{-2}+\mathbb{Z}\alpha^{-1}}\mathcal{R} \cap (\mathcal{R}+u)$.

On the other hand, if $z \in \mathcal{R} \cap (\mathcal{R}+u)$, $u \in \mathbb{Z}\alpha^{-2}+\mathbb{Z}\alpha^{-1}\backslash\{0\}$, since int$(\mathcal{R})\cap(\mathcal{R}+u)=\emptyset$, then $z \notin $ int($\mathcal{R}$). Then $z \in \partial \mathcal{R}$. Therefore, $\partial \mathcal{R}= \bigcup_{u \in \mathbb{Z}\alpha^{-2}+\mathbb{Z}\alpha^{-1}} \mathcal{R} \cap (\mathcal{R}+u)= \bigcup_{u \in H} \mathcal{R} \cap (\mathcal{R}+u)$, where $H=\{u \in \mathbb{Z}\alpha^{-2}+\mathbb{Z}\alpha^{-1}\backslash \{0\}, \mathcal{R} \cap  (\mathcal{R}+u)\neq \emptyset\}$.

The set $H$ is finite because $H$ is bounded and it belongs to  $\mathbb{Z}\alpha^{-2}+\mathbb{Z}\alpha^{-1}$, which is a lattice. Finally, the cardinality of $H$ is even because if $u \in H$, then $-u \in H$.

Let us prove that $\{\pm (1+(b+1)\alpha, \pm \alpha,\pm (1+b\alpha)\} \subset H$. By Lemma \ref{pontow}, we have seen that $w_{1} \in \mathcal{R}\cap (\mathcal{R}+\alpha) \cap (\mathcal{R}+1+(b+1)\alpha)$. Therefore $-\alpha$ and $-1-(b+1)\alpha$ belong to $H$. We have also seen that $z_{2} \in \mathcal{R} \cap (\mathcal{R}+1+b\alpha)$. Therefore, $-1-b\alpha$ belongs to $H$. $\Box$

\begin{remark} \label{triplo} We have seen in Lemma \ref{pontow} that a point, for instance $w_{1}$, belongs to $ \mathcal{R} \cap (\mathcal{R}+\alpha)$. This means that $w_{1}$ has two ways to be represented. Actually, in that case, $w_{1}$ could be expressed in three different ways. Points like $w_{1}$ are said to have at least two $\alpha$-representations. These points will be characterized in the next section. 
\end{remark}

\begin{figure}[h!]
\centering
 \includegraphics[scale=.34]{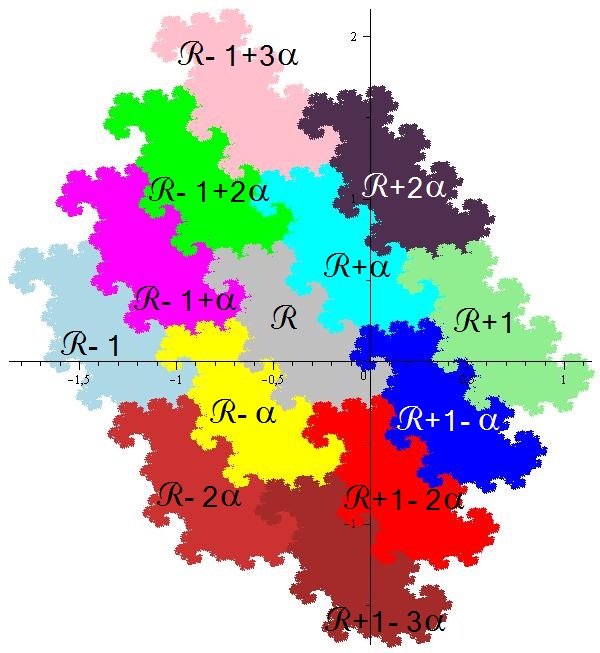}
 \includegraphics[scale=.29]{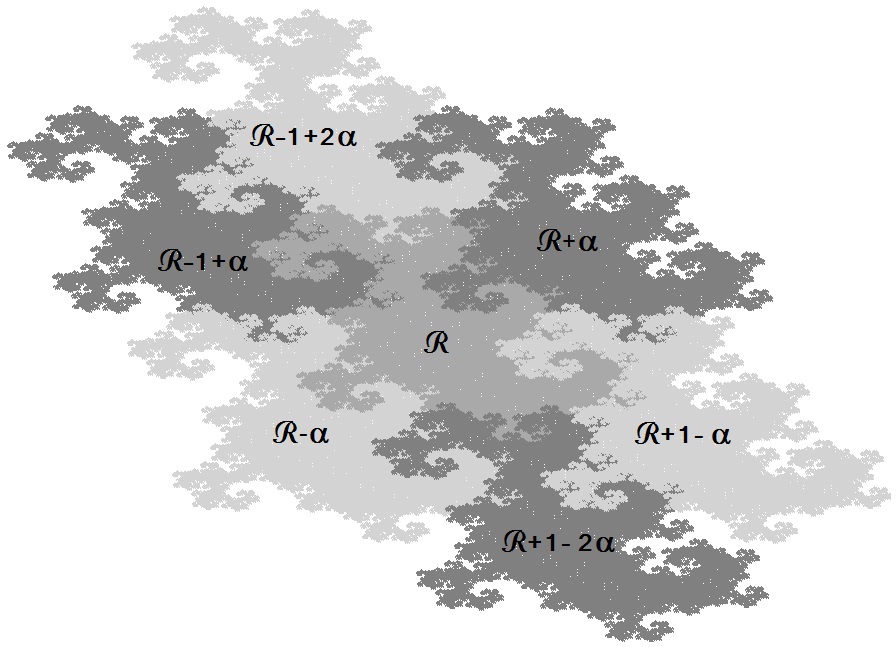}
\centerline{}
\centerline{Figure 3: Tiling the plane by $\mathcal{R}_{4,-3}$. \hspace{2.0cm} Figure 4: $\mathcal{R}_{3,-2}$ and its 6 neighbors}

\end{figure}

\section{Construction of the automaton $\cal{G}$}

In this section we prove that there exists an explicit and finite automata that recognize the points with two representations. These points belong to the boundary of $\mathcal{R}$. Let us begin with the following result.

\begin{prop} \label{xyteorema} Let $x=\sum_{i=l}^{\infty}a_{i}\alpha^{i}$ and $y=\sum_{i=l}^{\infty}b_{i}\alpha^{i}$, where  $l \in \mathbb{Z}$ and $(a_{i})_{i\geq l}$, $(b_{i})_{i\geq l}$ belong to $\mathcal{L}$. Then $x=y$ if, and only if, the set  $J(x,y)=\{x(k)-y(k), k \geq l\}$ is finite, where $x(k)=\alpha^{-k+2}\displaystyle \sum_{i=l}^{k}a_{i}\alpha^{i}$ and $y(k)=\alpha^{-k+2} \sum_{i=l}^{k}b_{i}\alpha^{i}$, $\forall k \geq l$.

Moreover, $\bigcup_{(x, y)}J(x,y) \supset E_{a,b}=\{0, \alpha^{2}, (\alpha+b\alpha^{2}), (\alpha +(b+1)\alpha^{2}), (1+b\alpha+(a-1)\alpha^{2}),  (1+(b+1)\alpha+(a+b)\alpha^{2}),  (1+(b+1)\alpha+(a+b+1)\alpha^{2})\}$.
\end{prop}

Before proving the Proposition, we will construct the automaton.  

\subsection{Algorithmic construction of the complex numbers that have two representations}

Let $p$ and $q$ be two states. The set of the edges is the set of $(p,(c,d),q) \in E_{a,b}\times\{0,1,...,a-1\}^{2}\times E_{a,b}$ satisfying $q=\frac{p}{\alpha}+(c-d)\alpha^{2}$. The set of the initial states is $\{0,(0,0),0\}$. 

Let us explain how this automaton acts. Let $x=\sum_{i=l}^{+\infty}a_{i}\alpha^{i}$ and $y=\sum_{i=l}^{+\infty}b_{i}\alpha^{i}$, where $a=(a_{i})_{i \geq l}$ and $b=(b_{i})_{i \geq l}$ belong to $\cal{L}$. Suppose that $x=y$ and for all $k \geq l$ we set $S_{k}=S_{k}(a,b)=x(k)-y(k)$. We have,

 \begin{equation} \label{Akmais1} \hspace{4.0cm }\displaystyle S_{k+1}=\frac{S_{k}}{\alpha}+(a_{k+1}-b_{k+1})\alpha^{2}. \end{equation}

Let $t$ be the smallest integer such that $a_{t} \neq b_{t}$. Hence $S_{i}(a,b)=0$ for all $i \in \{l,...,t-1\}$. Suppose that $a_{t}>b_{t}$. Then, $S_{t}=(a_{t}-b_{t})\alpha^{2}=\alpha^{2}$. From (\ref{Akmais1}) we deduce that $S_{t+1}=\alpha+(a_{t+1}-b_{t+1})\alpha^{2}$ which should belong to $E_{a,b}$. Hence $S_{t+1}=\alpha+b\alpha^{2}$ if $(a_{t+1},b_{t+1})=(s_{1}+b,s_{1})$, where $0 \leq s_{1} \leq a-1$, or $S_{t+1}=\alpha+(b+1)\alpha^{2}$, if $(a_{t+1},b_{t+1})=(t_{1}+b+1,t_{1})$, where $0 \leq t_{1} \leq a-1$. Continuing with this process, we obtain an infinite path $(S_{i}, (a_{i},b_{i}),S_{i+1})_{i \geq l}$ beginning in the initial state of the finite state automaton (see Fig. 5). This path will be denoted by $(a_{i},b_{i})_{i \geq l}$.\\

\begin{center}
   \includegraphics[width=12cm,angle=0]{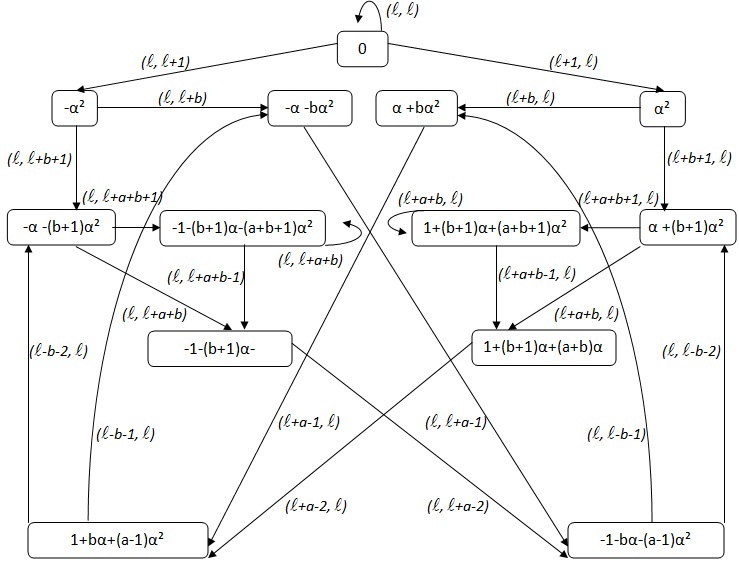}
\end{center} \centerline{Figure 5. Automaton $\mathcal{G}$}

\bigskip

\noindent \textbf{Proof of Proposition \ref{xyteorema}.} The direct implication is easy to see. Let us prove the converse. Let $x= \sum_{i=l}^{\infty}a_{i}\alpha^{i}$ and $y = \sum_{i=l}^{\infty}b_{i}\alpha^{i}$. Suppose that $x=y$, then $\alpha^{-k+2}x= \alpha^{-k+2}y$. Let us prove that the set $\{x(k)-y(k), k \geq 0\}$ is finite. Since $x(k)-y(k)=\alpha^{-k+2}(\sum_{i=0}^{k}a_{i}\alpha^{i}-\sum_{i=0}^{k}b_{i}\alpha^{i})=\alpha^{-k+2}(\sum_{i=k+1}^{+\infty}b_{i}\alpha^{i}-\sum_{i=k+1}^{+\infty}a_{i}\alpha^{i})=\sum_{j=3}^{+\infty}(b_{k+j-2}-a_{k+j-2})\alpha^{j}$, then $|x(k)-y(k)|\leq C$, where $C >0$ is a constant.

Let $S_{k} = x(k)-y(k)$. Then $S_{k}$ is an algebraic integer whose conjugates are $\widetilde{S}_{k}$ and $\overline{S}_{k}$, where $\widetilde{S}_{k}=\sum_{i=0}^{k}(a_{i}-b_{i})\beta^{i-k+2}$ and $\overline{S}_{k}=\sum_{i=0}^{k}(a_{i}-b_{i})\overline{\alpha}^{i-k+2}$.

We have $|\overline{S}_{k}|=|S_{k}| \leq C$, and \\

$|\widetilde{S}_{k}|=\left| \sum_{i=l}^{k}(a_{i}-b_{i})\beta^{i-k+2}\right|=|(a_{0}-b_{0})\beta^{-k+2}+ \cdots + (a_{k}-b_{k})\beta^{2}| \leq C \frac{\beta^{2}}{1-(1/\beta)}$,\\

where $C=2\cdot max\{|a_{i}|,a_{i} \in \{0,1,\ldots, a-1 \}\}=2(a-1)$.

Then, there exists $M>0$ such that $S_{k}$ and all its conjugates are bounded by $M$, independently of $k$. Thus $\{S_{k},k\geq 0\}$ is finite. $\Box$\\

As a consequence of this Proposition, we have the following result.

\begin{theo} Let $(a_{i})_{i \geq l}$ and $(b_{i})_{i \geq l}$ two distinct elements of $\mathcal{L}$, then $\sum_{i=l}^{\infty}a_{i}\alpha^{i}=\sum_{i=l}^{\infty}b_{i}\alpha^{i}$ if and only if the sequence $((a_{i},b_{i}))_{i \geq l}$ is recognizable by the automaton $\mathcal{G}$.
\end{theo}

\begin{remark} The usage of finite state automata to recognize points with two expansions is well known (see \cite{Thurston},\cite{Frougny2000},\cite{Me2000}). The difficult is to find the states of these automata. And that is what we are going to do in the sequel.
\end{remark}

To prove that $\bigcup_{(x, y)}J(x,y) \supset E_{a,b}$ we need the following result:

\begin{prop} \label{states} Let $F_{a,b}=\{S_{k}=n_{k}+p_{k}\alpha+q_{k}\alpha^{2}, k \geq 0$ and $(n_{k},p_{k},q_{k})\in I$, where $I \subset \mathbb{Z}\times\mathbb{Z}\times\mathbb{Z}$ is finite$\}$ and $t=max\{|n_{k}|,k \geq 0\}$. If $t=1$, then the set of the states of the automaton contains $E_{a,b}$.
\end{prop}

For proving Proposition \ref{states}, we need the next Lemma.

\begin{lemm} \label{betacubo} For all $k \geq l$, $|\widetilde{S}_{k}|=\left|\sum_{i=l}^{k}(a_{i}-b_{i})\beta^{i-k+2}\right|<\beta^{3}$.
\end{lemm}

\noindent \textbf{Proof.} Suppose, without loss of generality, that $\widetilde{S}_{k}= \sum_{i=l}^{k}(a_{i}-b_{i})\beta^{i-k+2}>0$. Then, $\widetilde{S}_{k} \in \mathbb{Z}[\beta] \cap \mathbb{R}^{+}$. Since $\beta$ is a Pisot number, $\widetilde{S}_{k}=  \sum_{i=-\infty}^{L}c_{i}\beta^{i}$, where $(c_{i})_{i \leq L}$ is ultimately periodic (see \cite{Klaus}). Then, $\sum_{i=l}^{k}a_{i}\beta^{i-k+2}=\sum_{i=l}^{k}b_{i}\beta^{i-k+2}+\sum_{-\infty}^{L}c_{i}\beta^{i}$. Now, let us suppose that there exists $i \geq 3$ such that $c_{i}>0$. Then $ \sum_{i=l}^{k}a_{i}\beta^{i-k+2} \geq \beta^{3}$. Absurd, because $0a_{2}\cdots a_{l} <_{lex} 10\cdots 0$. Hence $L \leq 2$ and $\widetilde{S}_{k}<\beta^{3}$. $\Box$\\

\noindent \textbf{Proof of Proposition \ref{states}.} Let $S_{k}=n_{k}+p_{k}\alpha+q_{k}\alpha^{2}$ and $t= max\{|n_{k}|,k \geq 0\}$. Let us suppose that $t=1$. Then there exists and integer $k$ such that $S_{k}=1+p\alpha+q\alpha^{2}$. Then, by \ref{Akmais1}, $S_{k+1}=\frac{1}{\alpha}+p+q\alpha^{2}$. Hence, $S_{k+1}=(p-b)+(q-a)\alpha+(d+1)\alpha^{2}$, where $d \in \Lambda=\{-a+1,...,a-1\}$. Since $t=1$, then $ p \in \{b-1, b, b+1\}$.

\noindent Now we have to analyze all the possible values for $p$. Let us recall that $\beta^3=a\beta^2+b\beta+1$. \\

\noindent \textbf{Case 1.} $p=b$. In this case, $S_{k+1}=(q-a)\alpha+(d+1)\alpha^{2}$ and $S_{k+2}=(q-a)+(d+1)\alpha +e\alpha^{2}$, $e \in \Lambda$. Then, $ q \in \{a-1, a, a+1\}$. We have $\tilde{S}_{k}=1+b\beta+q\beta^{2}$. By Lemma \ref{betacubo}, we must have that $\tilde{S}_{k}<\beta^{3}$, hence $q=a-1$ because, otherwise, $\widetilde{S}_{k}\geq 1+b\beta+q\beta^{2}=\beta^{3}$. Hence we have the state $S_{k}=1+b\alpha+(a-1)\alpha^{2}$.\\

\noindent \textbf{Case 2.} $p=b-1$. We have: $S_{k+1}=-1+(q-a)\alpha+(d+1)\alpha^2$. Since $-S_{k+1}=1-(q-a)\alpha-(d+1)\alpha^{2} \in F_{a,b}$, we obtain as before that $a-q\in \{b-1, b, b+1\}$. Let us show that these cases do not occur. \\

\indent \textbf{2.1.} $q=a-b+1$. In this case, $\tilde{S}_{k}=1+(b-1)\beta+(a-b+1)\beta^{2}=\beta^{3}-\beta+(1-b)\beta^{2} \geq \beta^{3}-\beta+3\beta^{2}>\beta^{3}$. Hence $\tilde{S}_{k}>\beta^{3}$.\\

\indent \textbf{2.2.} $q=a-b$. We have $\tilde{S}_{k}=1+(b-1)\beta+(a-b)\beta^{2}=\beta^{3}-\beta -b\beta^{2} \geq \beta^{3}-\beta+2\beta^{2} > \beta^{3}$.\\

\indent \textbf{2.3.} $q=a-b-1$. In this case, $\tilde{S}_{k}=1+(b-1)\beta+(a-b-1)\beta^{2}=\beta^{3}-\beta+(-1-b)\beta^{2} \geq \beta^{3}-\beta+\beta^{2}>\beta^{3}$. So we do not have the case when $p=b-1$.\\

\noindent \textbf{Case 3.} $p=b+1$. We have $S_{k+1}=1+(q-a)+(d+1)\alpha^{2}$ and $S_{k+2}=(q-a-b)+(d+1-a)\alpha+(e+1)\alpha^{2}$. Then, $q=a+b+r$, where $|r|\leq 1$. Hence, 

\centerline{$\tilde{S}_{k}=1+p\beta+q\beta^{2}=1+(b+1)\beta+(a+b+r)\beta^{2}$, }

with $|r|\leq 1$. We have to analyze all the possible cases for $q$.\\

\indent \textbf{3.1.} $q=a+b$. In this case, $\tilde{S}_{k}=1+(b+1)\beta+(a+b)\beta^{2}=\beta^{3}+\beta+b\beta^{2} \leq \beta^{3}+\beta-2\beta^{2}< \beta^{3}$. So we have the state $S_{k}=1+(b+1)\alpha+(a+b)\alpha^{2}$.\\

\indent \textbf{3.2.} $q=a+b+1$. We have: 

\centerline{$\tilde{S}_{k}=1+(b+1)\beta+(a+b+1)\beta^{2}=\beta^{3}+\beta+(b+1)\beta^{2} \leq \beta^{3}+\beta-\beta^{2}< \beta^{3}$. }

Then we obtain the state $S_{k}=1+(b+1)\alpha+(a+b+1)\alpha^{2}$.\\

\indent \textbf{3.3.} $q=a+b-1$. In this case,  $\tilde{S}_{k}=1+(b+1)\beta+(a+b-1)\beta^{2}.$ Hence $S_{k+1}=1+(b-1)\alpha+(d+1)\alpha^{2}$, $d \in \Lambda$, which does not exist by \textbf{Case 2}. So, this case does not occur. \\

\noindent Let us now consider $S_{k}=n+p\alpha+q\alpha^{2}$ and suppose that $n=0$. Then, $S_{k}=p\alpha+q\alpha^{2}$ and $A_{k+1}=p+q\alpha+d\alpha^{2}$. Then, $p \in \{-1,0,1\}$. Let us analyze all the possible cases, as we have done before.\\

\noindent \textbf{Case 4.} If $p=0$ then $S_{k}=q\alpha^{2}$. So we obtain the states $S_{k}=0$, if $q=0$, and $S_{k}=\pm\alpha^{2}$, if $q= \pm 1$.\\

\noindent \textbf{Case 5.} If $p=1$, then $S_{k}=\alpha+q\alpha^{2}$ and $S_{k+1}=1+q\alpha+d\alpha^{2}$. Hence, $S_{k+2}=(q-b)+(d-a)\alpha+(e+1)\alpha^{2}$. Thus, $q \in \{b-1,b,b+1 \}$. Let us analyze the possible cases. \\

\indent \textbf{5.1.} $q=b-1$. This case does not occur, as seen in \textbf{Case 2}.\\

\indent \textbf{5.2.} $q=b$. In this case we have the state $S_{k}=\alpha+b\alpha^{2}$.\\

\indent \textbf{5.3.} $q=b+1$. In this case we have the state $S_{k}=\alpha+(b+1)\alpha^{2}$. $\Box$\\

Next Proposition, which is a very important one, tell us that the automaton could have other states depending on certain conditions.\\

\begin{prop} \label{morestates} Let $t$ be the integer defined in Proposition \ref{states} and  suppose that  $1< t \leq \frac{a-1}{a+b+1}$. If $a+b \geq 3$ then $S_{k}=t+(t+tb)\alpha+(ta+tb+t)\alpha^{2}$ is a state of the automaton $\mathcal{G}$.
\end{prop}

\noindent \textbf{Proof.} This proof highly depends on the properties of the associated \linebreak $\beta-$expansion and it must be divided into several cases. Let us remind that $d(1,\beta)=(a-1)(a+b-1)(a+b)^{\infty}$.\\

\noindent \textbf{Case 1.} $t \geq 2$. Suppose that $S_{k}=t+p\alpha+q\alpha^{2}$. Then, by (\ref{Akmais1}), $S_{k+1}=(p-tb)+(q-ta)\alpha+(f+t)\alpha^{2}$, where $|f|\leq a-1$. Setting $p-tb=s$, where  $s \in \{-t,\ldots,t\}$, thus $S_{k+1}=t+(q-ta)\alpha+(f+t)\alpha^{2}$. Then, $S_{k+2}=(q-ta-sb)+(f+t-sa)\alpha+(g+s)\alpha^{2}$, where $|g|\leq a-1$. Setting $q=ta+sb+l$, for $l=-t,\ldots,t$ we obtain that $S_{k}=t+(tb+s)\alpha+(ta+sb+l)\alpha^{2}$. Let us show that $s=t$ whenever $a+b\geq 3$.\\

We have $\widetilde{S}_{k}=t+(tb+s)\beta+(ta+sb+l)\beta^{2}$. Since $t+tb\beta+ta\beta^{2}=t\beta^{3}$, we obtain that $\widetilde{S}_{k}=t\beta^{3}+(sb+l)\beta^{2}+s\beta=\beta^{3}+(t-1)\beta^{3}+(sb+l)\beta^{2}+s\beta$. 

Set $X=(t-1)\beta^{3}+(sb+l)\beta^{2}+s\beta$. Using the fact that $\beta^{2}=(a-1)\beta+(a+b-1)+(a+b)\sum_{i=1}^{\infty}1/\beta^{i}$, we obtain

$\begin{array}{ll}
  X/\beta & =(t-1)\beta^{2}+(sb+l)\beta+s \\
   & =[(t-1)(a-1)+sb+l]\beta+[(t-1)(a+b-1)+s] +R, \\
\end{array}$

where $R=(t-1)(a+b)\sum_{i=1}^{\infty}1/\beta^{i}>0$.\\

By Lemma \ref{betacubo} we must have $\widetilde{S}_{k}<\beta^{3}$. So, we need to show that $X/\beta \geq 0$. Let us do the first two cases. For all cases, the reader is referred to \cite{TeseGustavo}).

\noindent Let us suppose that $s < t$.\\

\noindent \textbf{Case 1.1.} $-t \leq s \leq 0$. In this case,

$\begin{array}{ll}
  X/\beta & =[(t-1)(a-1)+sb+l]\beta+[(t-1)(a+b-1)+s]+R \\
   & \geq [l+(t-1)(a-1)]\beta+[(t-1)(a+b-1)+s]+R \hbox{, since } sb\geq 0 \\
   & \geq [(t-1)(a-1)-t]\beta+ [(t-1)(a+b-1)-t]+R \\
\end{array}$\\
 since $l \geq -t$ and $s\geq -t$.\\

\textbf{Case} \textbf{1.1.1.} $a+b-1 \geq 2$. Since $t \geq 2$, $a \geq 3$ and $R > 0$ we obtain that $X/\beta \geq 0$ and then $X\geq 0$. So $\widetilde{S}_{k}=\beta^{3}+X > \beta^{3}$, which is an absurd.\\

After analyzing all the possibles cases we conclude that $s=t$ and $l=t$.\\

\textbf{Case 2}. $1<m<t$. Let $S_{k}=m+p\alpha+q\alpha^{2}$, with $m<t$, then $S_{k+1}=(p-mb)+(q-ma)\alpha+(d+m)\alpha^{2}$. Hence, $S_{k+1}=s+(q-ma)\alpha+(d+m)\alpha^{2}$ and $S_{k+2}=(q-ma-sb)+(d+m-sa)\alpha+(g+s)\alpha^{2}$, $|s| \leq t$, $|q| \leq t$, $|l| \leq t$. So, $S_{k}=m+(mb+s)\alpha+(ma+sb+l)\alpha^{2}$. \\

\textbf{Remark.} Since $ta+tb+t \in \mathcal{L}$, then it must satisfies the condition: $0 \leq ta+tb+t \leq a-1$, that is, $\displaystyle n \leq \frac{a-1}{a+b+1}$. \\

Therefore, $S_{k}=t+(t+tb)\alpha+(ta+tb+t)\alpha^{2}$ is a state of the automaton $\mathcal{G}$. $\Box$\\

\begin{corol} The automaton has at least $2(6+2(K-1))$ nonempty states, where $\displaystyle K=\left[\frac{a-1}{a+b+1}\right]$. The set of the states contains $E_{a,b} \cup \{\pm t\alpha^{2}, \pm (t\alpha+t(b+1)\alpha^{2},\pm(t+t(b+1)\alpha+t(a+b+1)\alpha^{2})\}$, where $t \leq \displaystyle \left[ \frac{a-1}{a+b+1}\right]$.
\end{corol}

\begin{corol} \label{voisins} $\mathcal{R}$ has at least $6+2(K-1)$ neighbors of the form $u+\mathcal{R}$, where $\pm u \in \{\pm \alpha,\pm (1+b\alpha), \pm (1+(b+1)\alpha)\} \cup \{\pm (t+t(b+1)\alpha\}$ and $\displaystyle t \leq \left[ \frac{a-1}{a+b+1} \right]$.
\end{corol}

From Corollary \ref{voisins} we have the following Theorem. 

\begin{theo} If $2a+3b+4 \leq 0$ then $\mathcal{R}$ is not homeomorphic to a topological disk.
\end{theo}

\noindent \textbf{Proof.} If $2a+3b+4 \leq 0$ then $3 \leq \displaystyle \frac{a-1}{a+b+1}$, that is, $\displaystyle K=\left[\frac{a-1}{a+b+1}\right] \geq 3$. Thus $\mathcal{R}$ has at least $6+2(K-1) \geq 10$ neighbors. So $\mathcal{R}$ cannot be homeomorphic to a topological disk (see \cite{bandt}). $\Box$

\begin{center}
    \includegraphics[width=7cm,angle=0]{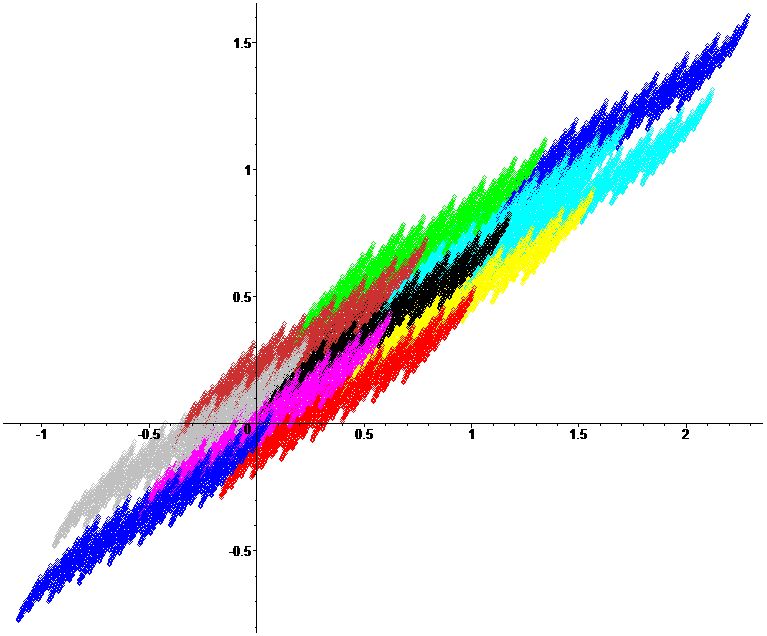}
  \end{center}
\centerline{Figure 6: $\mathcal{R}_{8,-7}$ has 10 neighbors.}

\section{Parametrization of the Boundary of $\mathcal{R}_{3,-2}$}

In this section we will use the Automaton $\mathcal{G}$ built in the previous section with $a=3$ and $b=-2$ to generate the boundary of $\mathcal{R}=\mathcal{R}_{3,-2}$. By Corollary \ref{voisins}, $\mathcal{R}_{3,-2}$ has only 6 neighbors. We will prove that the boundary of $\mathcal{R}$ is generated  by two infinite countable sets of IFS.\\

Let $u \in \{\pm\alpha,\pm(1-\alpha),\pm(1-2\alpha)\}$ and denote by $\mathcal{R}_{u}=\mathcal{R} \cap (\mathcal{R}+u)$ the 6 curves which constitute the boundary of $\mathcal{R}$. The next proposition shows that each neighbor of $\mathcal{R}$ can be expressed by means of the other ones.  

\begin{prop} \label{vizinhos} The following relations are valid:

1. $\mathcal{R}_{1-\alpha}=\bigcup_{k=1}^{\infty}\ell_{k+1}\alpha^{k+1}+\alpha^{k}\mathcal{R}_{1-2\alpha}$, where $\ell_{k+1} \in \{0,1,2\}$, for all $k \geq 0$.

2. $\mathcal{R}_{\alpha}=\alpha\mathcal{R}_{1-2\alpha}\bigcup \bigcup_{k=1}^{\infty}(\ell_{k+1}\alpha^{k+2}+\alpha^{k+1}\mathcal{R}_{1-2\alpha})$, where $\ell_{k+1} \in \{0,1,2\}$, for all $k \geq 1$.

\end{prop}

\noindent \textbf{Proof.} \textbf{1.} Let $z \in \mathcal{R}_{1-\alpha}$. Then $z=1-\alpha+\sum_{i=2}\ell_{i}\alpha^{i}=\sum_{i=2}d'_{i}\alpha^{i}$. So, by the Automaton  $\mathcal{G}$, we have the associated paths in the automaton beginning in the initial state: $P_{1}=(1,0)(-1,0)(\ell_{2}+1,\ell_{2})(1,0)\cdots$ or  

$P_{2}=(1,0)(-1,0)(2,0)\underbrace{(1,0)(1,0)\cdots(1,0)}_{k-times}(\ell_{3+k},\ell_{3+k})(1,0)\cdots$.\\

\noindent \textbf{Case 1.1} $z=1-\alpha+(\ell_{2}+1)\alpha^{2}+\alpha^{3}+\alpha^{4}w=\ell_{2}\alpha^{2}+\alpha^{4}w^{'}$, where $w,w^{'} \in \mathbb{C}$. Hence, $z/\alpha-\ell_{2}\alpha=1-2\alpha+2\alpha^{2}+\alpha^{3}w=\alpha^{3}w^{'} \in \mathcal{R}_{1-2\alpha}$.\\

On the other hand, if $z \in \mathcal{R}_{1-2\alpha}$ then $z=1-2\alpha+2\alpha^{2}+\alpha^{3}w$, where $w \in \mathbb{C}$. Thus, $\alpha z +\ell_{2}\alpha^{2}=\alpha+(\ell_{2}\alpha^{2}+2\alpha^{3}+\alpha^{4}w=1+(\alpha-2\alpha)+(\ell_{2}-3)\alpha^{2}+\alpha^{3}+\alpha^{4}w=1-\alpha=(\ell_{2}+1)\alpha^{2}+\alpha^{3}+\alpha^{4}w \in \mathcal{R}_{1-2\alpha}$, if $\ell_{2} \geq 1$.\\

\noindent \textbf{Case 1.2} $z=1-\alpha+2\alpha^{2}+\underbrace{\alpha^{3}+\cdots+\alpha^{k+2}}_{k-times}+\ell_{3+k}\alpha^{3+k}+\alpha^{4+k}+\alpha^{5+k}w_{k}=\ell_{3+k}\alpha^{3+k}+\alpha^{5+k}w_{k}^{'}$, where $w_{k},w_{k}^{'}\in\mathbb{C}$, for all $k \geq 0$. Hence, $z/\alpha^{k+2}=\alpha^{-k-2}-\alpha^{-k-1}+2\alpha^{-k}+\alpha^{-k+1}+\cdots+1+\ell_{3+k}\alpha+\alpha^{2}+\alpha^{3}w_{k}=\ell_{3+k}\alpha+\alpha^{3}w_{k}^{'}$. Thus, by induction, we can show that 

\centerline{$\displaystyle \frac{z}{\alpha^{k+2}}-\ell_{k+3}\alpha \in \mathcal{R}_{1-2\alpha}$.}

\bigskip

Therefore, $\mathcal{R}_{1-\alpha}=\bigcup_{k=1}^{\infty}\ell_{k+1}\alpha^{k+1}+\alpha^{k}\mathcal{R}_{1-2\alpha}$.\\

\noindent \textbf{2.} If $z\in \mathcal{R}_{\alpha}$ then $z=\alpha+\sum_{i=2}\ell_{i}\alpha^{i}=\sum_{i=2}\ell'_{i}\alpha^{i}$. Thus, $P_{1}=(0,0)(1,0)(-2,0)(2,0)\cdots$ or $P_{2}=(0,0)(1,0)(-1,0)\cdots$ are paths in the automaton beginning in the initial state.\\

\noindent \textbf{Case 2.1} $z=\alpha-2\alpha^{2}+2\alpha^{3}+\alpha^{4}w_{2}=\alpha^{4}w_{2}^{'}$. Hence, $z/\alpha=1-2\alpha+2\alpha^{2}+\alpha^{3}w_{2}=\alpha^{3}w_{2}^{'} \in \mathcal{R}_{1-2\alpha}$.\\

\noindent \textbf{Case 2.2} Then $z=\alpha-\alpha^{2}+\alpha^{3}w_{3}=\alpha^{3}w_{3}^{'}$. Hence, $z/\alpha=1-\alpha+\alpha^{2}w_{3}=\alpha^{2}w_{3}^{'} \in \mathcal{R}_{1-\alpha}$. (We are back in \textbf{Case 1.1}).\\

Therefore, $\mathcal{R}_{\alpha}=\alpha\mathcal{R}_{1-2\alpha} \cup \alpha\mathcal{R}_{1-\alpha}=\alpha\mathcal{R}_{1-2\alpha}\bigcup_{k=1}^{\infty}\ell_{k+1}\alpha^{k+2}+\alpha^{k+1}\mathcal{R}_{1-2\alpha}$. $\Box$\\

\noindent For all $z\in \mathbb{C}$, consider the following functions: 

$f_{0}(z)=-\alpha^{2}+2\alpha^{3}+\alpha^{2}z$; 

$f_{1,i,j}(z)=(i-1)\alpha^{3}+(j+3)\alpha^{4}-\alpha^{5}+\alpha^{3}z$, where $j=0$, if $i=1$, or $j\in \{0,1\}$, if $i=0$; 

$f_{2,i,j}(z)=i\alpha^{3}+\alpha^{4}+(j+2)\alpha^{5}-2\alpha^{6}+\alpha^{4}z$, where $i, j \in \{0,1\}$; 

$f_{3,i,j}(z)=i\alpha^{3}+2\alpha^{4}+(j+2)\alpha^{6}-2\alpha^{7}+\alpha^{5}z$, where $i, j \in \{0,1\}$;

$f_{2+k,i}(z)=2\alpha^{4}+(\sum_{j=1}^{k-1}\alpha^{4+j})+(i+2)\alpha^{5+k}-\alpha^{6+k}+\alpha^{4+k}z$, where $i\in \{0,1\}$ for all $k \geq 2$.\\

\begin{remark} Notice that some functions depend of the parameters $i$ and $j$. These functions will be separated into two sets (see Proposition \ref{y0z0}).  \end{remark}

The next Theorem shows that $\mathcal{R}_{1-2\alpha}$ is the infinite union of the images of itself by the applications defined above.  \\ 

\begin{theo} \label{r12alpha} $\mathcal{R}_{1-2\alpha}=\bigcup_{k=0}^{\infty}f_{k,i_{k},j_{k}}(\mathcal{R}_{1-2\alpha})$.

\end{theo}

\noindent \textbf{Proof.} For the first part of the proof, we need to show that each $f_{i} \subset \mathcal{R}_{1-2\alpha}$, for $i=0,1,...$. Let us do the computation for $f_{0}$, the other inclusions can be done in the same fashion. We have,

$f_{0}(\mathcal{R}_{1-2\alpha})=f_{0}(\mathcal{R}) \cap f_{0}(\mathcal{R}+1-2\alpha)=(-\alpha^{2}+2\alpha^{3}+\alpha^{2}\mathcal{R}) \cap (\alpha^{2}\mathcal{R})=(1-2\alpha+2\alpha^{2}+\alpha^{3}+\alpha^{2}\mathcal{R}) \cap (\alpha^{2}\mathcal{R}) \subset \mathcal{R}_{1-2\alpha}$.\\

On the other hand, let  $z \in \mathcal{R}_{1-2\alpha}$. Using the automaton $\mathcal{G}$ we have the following paths beginning in the initial state:

\noindent \textbf{1.} $P_{1}=(0,1)(0,-2)(0,2)(0,1)(2,0)\cdots$. Then, $z=2\alpha^{4}+\alpha^{5}w_{0}=1-2\alpha+2\alpha^{2}+\alpha^{3}+\alpha^{5}w_{0}^{'}$, where $w_{0},w_{0}^{'} \in \mathbb{C}$. Hence, $f_{0}^{-1}(z)=1-2\alpha+2\alpha^{2}+\alpha^{3}w_{0}=\alpha^{3}w_{0}^{'} \in \mathcal{R}_{1-2\alpha}$, that is, $z \in f_{0}(\mathcal{R}_{1-2\alpha})$.\\

\noindent \textbf{2.} $P_{2}=(0,1)(0,-2)(0,2)(\ell_{3},\ell_{3})(\ell_{4}+1,\ell_{4})(1,0)\cdots$. Then, $z=\ell_{3}\alpha^{3}+(\ell_{4}+1)\alpha^{4}+\alpha^{5}+\alpha^{6}w_{1}=1-2\alpha+2\alpha^{2}+\ell_{3}\alpha^{3}+\ell_{4}\alpha^{4}+\alpha^{6}w_{1}^{'}$. Thus, $f_{1}^{-1}(z)=1-2\alpha+2\alpha^{2}+\alpha^{3}w_{1}=\alpha^{3}w_{1}^{'} \in \mathcal{R}_{1-2\alpha}$, that is, $z \in f_{1}(\mathcal{R}_{1-2\alpha})$.\\

\noindent \textbf{3.} $P_{3}=(0,1)(0,-2)(0,2)(\ell_{3},\ell_{3})(2,0)(\ell_{5},\ell_{5})(1,0)\cdots$. Then, $z=\ell_{3}\alpha^{3}+2\alpha^{4}+\ell_{5}\alpha^{5}+\alpha^{6}+\alpha^{7}w_{2}=1-2\alpha+2\alpha^{2}+\ell_{3}\alpha^{3}+\ell_{5}\alpha^{5}+\alpha^{7}w_{2}^{'}$. Hence, $f_{2}^{-1}(z)=1-2\alpha+2\alpha^{2}+\alpha^{3}w_{2}=\alpha^{3}w_{2}^{'} \in \mathcal{R}_{1-2\alpha}$, that is, $z \in f_{2}(\mathcal{R}_{1-2\alpha})$.\\

\noindent \textbf{4.} $P_{4}=(0,1)(0,-2)(0,2)(\ell_{3},\ell_{3})(2,0)\underbrace{(1,0)(1,0)\cdots(1,0)}_{k-times}(\ell_{5+k},\ell_{5+k})(1,0)\cdots$. In this case, $z=\ell_{3}\alpha^{3}+2\alpha^{4}+\alpha^{5}+\alpha^{6}+\cdots+\alpha^{4+k}+\ell_{5+k}\alpha^{5+k}+\alpha^{6+k}+\alpha^{7+k}w_{2+k}=1-2\alpha+2\alpha^{2}+\ell_{3}\alpha^{3}+\ell_{5+k}\alpha^{5+k}+\alpha^{7+k}w_{2+k}^{'}$. Hence, $f_{2+k}^{-1}(z)=1-2\alpha+2\alpha^{2}+\alpha^{3}w_{2+k}=\alpha^{3}w_{2+k}^{'} \in \mathcal{R}_{1-2\alpha}$, for all $k \geq 1$, that is, $z \in f_{2+k}(\mathcal{R}_{1-2\alpha})$ for all $k\geq 1$.\\

Therefore $\mathcal{R}_{1-2\alpha}=\bigcup_{k=0}^{\infty}f_{k,i_{k},j_{k}}(\mathcal{R}_{1-2\alpha})$. $\Box$\\

We have shown that $\mathcal{R}_{1-2\alpha}= \bigcup_{a_{n},n \in \mathbb{N}}^{+\infty}f_{a_{n}}(\mathcal{R}_{1-2\alpha})$, where $a_{n}\in\{0,1,(2+k); k \in \mathbb{N}\}$. Now, let $z \in \mathcal{R}_{1-2\alpha}$. Then, $z=f_{a_{0}}(z_{0})=f_{a_{0}}\circ f_{a_{1}}\circ \cdots \circ f_{a_{n}}(z_{n})=\lim_{n \rightarrow +\infty}f_{a_{0}}\circ \cdots \circ f_{a_{n}}(z_{n})$, $z_{n}\in \mathcal{R}_{1-2\alpha}$ and $z$ is fixed. Thus, 

\centerline{$\displaystyle \mathcal{R}_{1-2\alpha}=\overline{\bigcup_{a_{0},\cdots, a_{n}} f_{a_{0}}\circ \cdots f_{a_{n}}(z)}.$}
\bigskip

Hence, using the Proposition \ref{vizinhos} and the Theorem \ref{r12alpha}, we obtain the boundary of $\mathcal{R}_{3,-2}$ (see Figure 11).\\

\noindent \textbf{Parametrization of $\mathcal{R}_{1-2\alpha}$}

The next Lemma shows points of the fractal $\mathcal{R}$ that can be expressed in three different ways. i.e.  a point with three $\alpha-$representations. Consequently these points lie in the intersection of three neighbors of the fractal (see Remark \ref{triplo}). These points are shown in Figure 7. 

\begin{lemm} The following properties are satisfied.

1. $\mathcal{R}_{\alpha} \cap \mathcal{R}_{1-\alpha} = \alpha + \frac{\alpha^{3}+\alpha^{4}+\alpha^{5}}{1-\alpha^{6}} $

2. $\mathcal{R}_{1-\alpha} \cap \mathcal{R}_{1-2\alpha} = \frac{\alpha^{4}+\alpha^{5}+\alpha^{6}}{1-\alpha^{6}}$.

3. $\mathcal{R}_{1-2\alpha} \cap \mathcal{R}_{-\alpha} = \frac{\alpha^{3}+\alpha^{4}+\alpha^{5}}{1-\alpha^{6}}$.

4. $\mathcal{R}_{-\alpha} \cap \mathcal{R}_{-1+\alpha} = -1+\alpha+\frac{\alpha^{4}+\alpha^{5}+\alpha^{6}}{1-\alpha^{6}}$.

5. $\mathcal{R}_{-1+\alpha} \cap \mathcal{R}_{-1+2\alpha}=-1+2\alpha + \frac{\alpha^{3}+\alpha^{4}+\alpha^{5}}{1-\alpha^{6}} $.

6. $\mathcal{R}_{-1+2\alpha} \cap \mathcal{R}_{\alpha} =-1+2\alpha+\frac{\alpha^{4}+\alpha^{5}+\alpha^{6}}{1-\alpha^{6}} $.

\end{lemm}

\textbf{Proof.} 1. If $w \in \mathcal{R}_{\alpha} \cap \mathcal{R}_{1-\alpha}$ then $w=\alpha+\sum_{i=2}^{+\infty}\ell_{i}\alpha^{i} = 1-\alpha+\sum_{i=2}^{+\infty}\ell^{'}_{i}\alpha^{i}$, where $(\ell_{i})_{i \geq 2},(\ell_{i}^{'})_{i \geq 2}\in \mathcal{L}$. Hence, using the automaton $\mathcal{G}$, we obtain that \\

$\begin{array}{ll}
  w & = 1-\alpha+2\alpha^{2}+\alpha^{3}+ \sum_{i=1}^{\infty}(\alpha^{6i+1}+\alpha^{6i+2}+\alpha^{6i+3})=\alpha^{2}+ \sum_{i=1}^{\infty}(\alpha^{6i-1}+\alpha^{6i}+\alpha^{6i+1})\\
   & =\alpha+ \sum_{i=1}^{\infty}(\alpha^{6i-3}+\alpha^{6i-2}+\alpha^{6i-1})=\alpha + \frac{\alpha^{3}+\alpha^{4}+\alpha^{5}}{1-\alpha^{6}}. \\
	\end{array}$ \\

2. If $x \in \mathcal{R}_{1-\alpha} \cap \mathcal{R}_{1-2\alpha}$ then $x=1-\alpha+\sum_{i=2}^{+\infty}\ell_{i}\alpha^{i} = 1-2\alpha+\sum_{i=2}^{+\infty}\ell^{'}_{i}\alpha^{i}$, where $(\ell_{i})_{i \geq 2},(\ell_{i}^{'})_{i \geq 2}\in \mathcal{L}$. Using the Automaton we obtain that\\

$\begin{array}{ll}
  x & =1-\alpha+ \sum_{i=1}^{\infty}(\alpha^{6i-4}+\alpha^{6i-3}+\alpha^{6i-2})=1-2\alpha+2\alpha^{2}+ \sum_{i=1}^{\infty}(\alpha^{6i-4}+\alpha^{6i-3}+\alpha^{6i-2})\\
   & =\sum_{i=1}^{\infty}(\alpha^{6i-2}+\alpha^{6i-1}+\alpha^{6i})=\frac{\alpha^{4}+\alpha^{5}+\alpha^{6}}{1-\alpha^{6}}. \\
	\end{array}$ \\

The other relations come from the fact that: $\mathcal{R}_{1-2\alpha} \cap \mathcal{R}_{-\alpha} = \mathcal{R}_{\alpha} \cap \mathcal{R}_{1-\alpha} - \alpha$, $\mathcal{R}_{-\alpha} \cap \mathcal{R}_{-1+\alpha} = \mathcal{R}_{1-\alpha} \cap \mathcal{R}_{1-2\alpha} -1+\alpha$, $\mathcal{R}_{-1+\alpha} \cap \mathcal{R}_{-1+2\alpha}= \mathcal{R}_{\alpha} \cap \mathcal{R}_{1-\alpha}-1+\alpha$, and $\mathcal{R}_{-1+2\alpha} \cap \mathcal{R}_{\alpha} = \mathcal{R}_{1-\alpha} \cap \mathcal{R}_{1-2\alpha} -1 +2\alpha$.$\Box$

\begin{center}
    \includegraphics[width=12cm,angle=0]{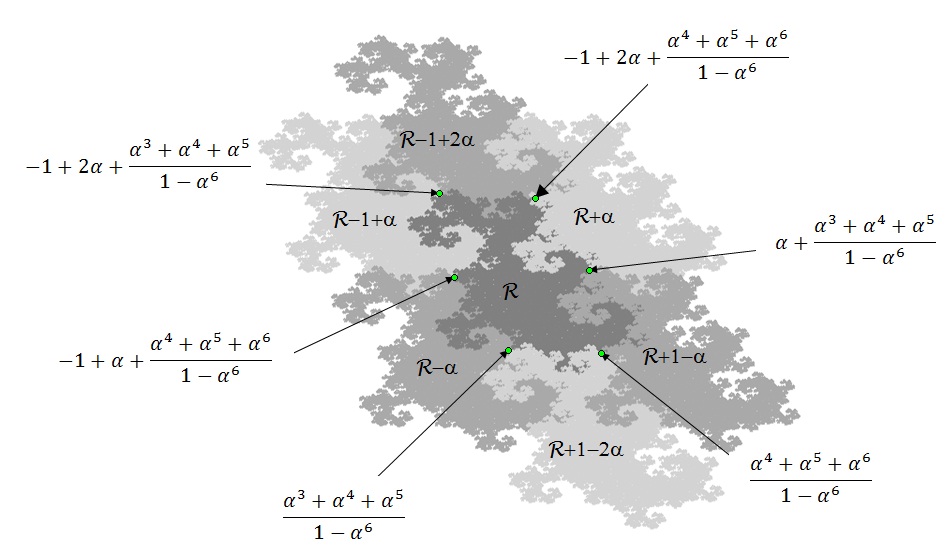}
  \end{center}
\centerline{Figure 7: Points with three $\alpha$-representations.}

\bigskip

Let us consider the iterated function system consisting of:

$f_{0}(z)=3\alpha^{4}-\alpha^{5}+\alpha^{3}z$, 

$f_{1}(z)=-\alpha^{2}+2\alpha^{3}+\alpha^{2}z$, 

$f_{2}(z)=-\alpha^{3}+4\alpha^4-\alpha^{5}+\alpha^{3}(z)$, 

$f_{3}(z)=\alpha^{4}+3\alpha^{5}-\alpha^{6}+\alpha^{4}(z)$, 

$f_{4}(z)=2\alpha^{4}+3\alpha^{6}-\alpha^{7}+\alpha^{5}(z)$, 

$f_{3+k}(z)=2\alpha^{4}+(\sum_{j=1}^{k-1}\alpha^{4+j})+3\alpha^{5+k}-\alpha^{6+k}+\alpha^{4+k}z$, for all $k \geq 2$;\\

and \\

$g_{0}(z)=-\alpha^{3}+3\alpha^{4}-\alpha^{5}+\alpha^{3}z$, 

$g_{1}(z)=\alpha^{4}+2\alpha^{5}-\alpha^{6}+\alpha^{4}z$, 

$g_{2}(z)=2\alpha^{4}+2\alpha^{6}-\alpha^{7}+\alpha^{5}z$, 

$g_{1+k}(z)=2\alpha^{4}+(\sum_{j=1}^{k-1}\alpha^{4+j})+2\alpha^{5+k}-\alpha^{6+k}+\alpha^{4+k}z$, for all $k \geq 2$.\\

Figure 8 illustrates the behaviour of this system.  We have the following result.

\begin{prop} \label{y0z0} For all $i,l \in \mathbb{N}$,

1. $f_{i}(\mathcal{R}_{1-2\alpha}) \cap f_{l}(\mathcal{R}_{1-2\alpha})\neq \emptyset$ if, and only if, $0 \leq |i-l| \leq 1$. In particular, $f_{3+k}(\mathcal{R}_{1-2\alpha}) \cap f_{3+(k+1)}(\mathcal{R}_{1-2\alpha})=\{f_{3+k}(z_{0})\}=\{f_{3+(k+1)}(y_{0})\}$, where $z_{0}=\frac{\alpha^{3}+\alpha^{4}+\alpha^{5}}{1+\alpha^{6}}$ and $y_{0}=\frac{\alpha^{4}+\alpha^{5}+\alpha^{6}}{1-\alpha^{6}}$;

\bigskip

2. $ g_{i}(\mathcal{R}_{1-2\alpha}) \cap g_{l}(\mathcal{R}_{1-2\alpha})\neq \emptyset$ if, and only if, $0 \leq |i-l| \leq 1$. In particular, $g_{1+k}(\mathcal{R}_{1-2\alpha}) \cap g_{1+(k+1)}(\mathcal{R}_{1-2\alpha})=\{g_{1+k}(y_{0})\}=\{g_{1+(k+1)}(z_{0})\}$, where $z_{0}=\frac{\alpha^{3}+\alpha^{4}+\alpha^{5}}{1+\alpha^{6}}$ and $y_{0}=\frac{\alpha^{4}+\alpha^{5}+\alpha^{6}}{1-\alpha^{6}}$;

\bigskip

3.  $f_{i}(\mathcal{R}_{1-2\alpha}) \cap g_{l}(\mathcal{R}_{1-2\alpha}) = \emptyset$, for all $i,l \in \mathbb{N}$.   

\end{prop}

\textbf{Proof.} Let us prove the item \textit{1}. \\

\mathversion{bold} \textbf{Case:} $0 \leq |i-l| \leq 1$. \mathversion{normal}

Let us suppose that $w \in f_{3+k}(\mathcal{R}_{1-2\alpha}) \cap f_{3+(k+1)}(\mathcal{R}_{1-2\alpha})$. Then there exists $y,z \in \mathcal{R}_{1-2\alpha}$ such that $y=-\alpha+\alpha^{2}+\alpha z \in \mathcal{R}_{-\alpha} \cap \mathcal{R}_{1-2\alpha}$. Hence, $y=\{z_{0}\}$ and $z=\{y_{0}\}$. Therefore, $f_{3+k}(\mathcal{R}_{1-2\alpha}) \cap f_{3+(k+1)}(\mathcal{R}_{1-2\alpha})=\{f_{3+k}(z_{0})\}=\{f_{3+(k+1)}(y_{0})\}$. 

In the same way, we can show that $f_{0}(\mathcal{R}_{1-2\alpha}) \cap f_{1}(\mathcal{R}_{1-2\alpha})=\{f_{0}(z_{0})\}=\{f_{1}(y_{0})\}$, $f_{1}(\mathcal{R}_{1-2\alpha}) \cap f_{2}(\mathcal{R}_{1-2\alpha})=\{f_{1}(z_{0})\}=\{f_{2}(y_{0})\}$, and $f_{2}(\mathcal{R}_{1-2\alpha}) \cap f_{3}(\mathcal{R}_{1-2\alpha})=\{f_{2}(z_{0})\}=\{f_{3}(y_{0})\}$.\\

\mathversion{bold}\textbf{Case:} $|i-l|>1$.\mathversion{normal}

Suppose that $l >i$ and that $f_{i}(\mathcal{R}_{1-2\alpha}) \cap f_{l}(\mathcal{R}_{1-2\alpha}) \neq \emptyset$. Then there exists $y,z \in \mathcal{R}_{1-2\alpha}$ such that 
\begin{equation} \label{intersect}
\sum_{j=1}^{i-1}\alpha^{4+j}+3\alpha^{5+i}-\alpha^{6+i}+\alpha^{4+i}y=\sum_{j=1}^{l-1}\alpha^{4+j}+3\alpha^{5+l}-\alpha^{6+l}+\alpha^{4+l}z.
\end{equation}

Since $y,z \in \mathcal{R}_{1-2\alpha}$, they can be expressed as $y=1-2\alpha+2\alpha^{2}+\alpha^{3}\bar{y}$ and $z=1-2\alpha+2\alpha^{2}+\alpha^{3}\bar{z}$, where $\bar{y},\bar{z} \in \mathbb{C}$. Replacing this in the equation (\ref{intersect}) we obtain that 

\begin{equation} \label{intersect2}
\bar{y}=1+\alpha+\alpha^{2}+ \cdots +\alpha^{l-i-1}+\alpha^{l-i}(\bar{z}). 
\end{equation}

Thus, 

\centerline{$\underbrace{(1,0)(1,0)(1,0)\ldots(1,0)}_{(l-i-1)\,\,times}\ldots$}

is the associated path in the automaton begining in the initial state that represents the point in (\ref{intersect2}). Absurd, because there is no such a path in the automaton. 

Therefore, $f_{i}(\mathcal{R}_{1-2\alpha}) \cap f_{l}(\mathcal{R}_{1-2\alpha}) = \emptyset$.\\

Using the same reasoning we can prove the items \textit{2.} and \textit{3.} $\Box$

\begin{center}
    \includegraphics[scale=.30]{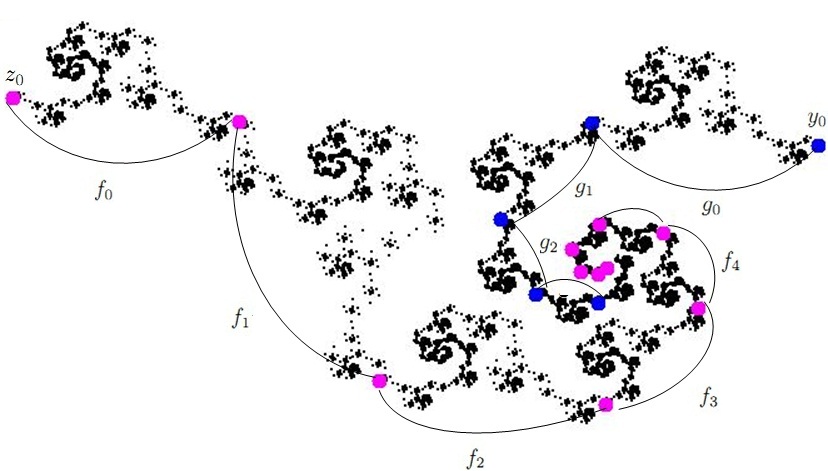}
  \end{center}
\centerline{Figure 8: $\mathcal{R}_{1-2\alpha}$}

\bigskip

Now we show a geometric way for parametrizing $\mathcal{R}_{1-2\alpha}$.  Figure 9 illustrates this procedure. Let $z_{0}$ and $y_{0}$ be two end points of $\mathcal{R}_{1-2\alpha}$ as in proposition \ref{y0z0}. Let us consider the sequence of function $\varphi_{n}:[0,1] \rightarrow \mathbb{C}, n \geq 1$, where:\\

$\varphi_{1}([0,1])$  is the polygonal line made up of segments of the form $[f_{k}(y_{0}),\, f_{k+1}(y_{0})]$, for $k \in \mathbb{N}$, and the segments $[g_{k}(z_{0}),\, g_{k+1}(z_{0})]$, for $k \in \mathbb{N}$. Let us remark that they could be joint in a continuous way (see last Proposition). \\

$\varphi_{2}([0,1])$ is the polygonal line consisting of all of the segments $[f_{i} \circ f_{j}(x), f_{i} \circ f_{j'}(x)]$, $[f_{i} \circ g_{j}(x), f_{i} \circ g_{j'}(x)]$, where $x\in\{z_{0},y_{0}\}$, and $i,j,j' \in \mathbb{N}$. See Figure 10 for clarity. \\

Once $\varphi_{n}([0,1])$ has been constructed, since each $f_{i}$ and $g_{j}$ are contractions, it is possible to show that $\varphi_{n}:[0,1] \rightarrow \mathbb{C}$ converges uniformly to a continuous function $h:[0,1] \rightarrow \mathcal{R}_{1-2\alpha}$. \\

Notice that with this method we can parametrize the whole boundary of $\mathcal{R}_{3,-2}$, once each neighbor is expressed by means of the others. \\

\begin{figure}[h]
\centering
{} \includegraphics[scale=0.35]{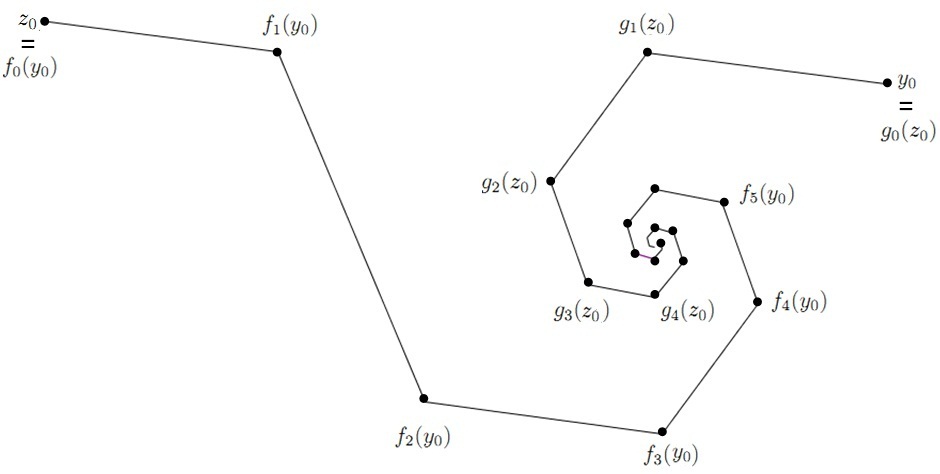}
{} \includegraphics[scale=0.35]{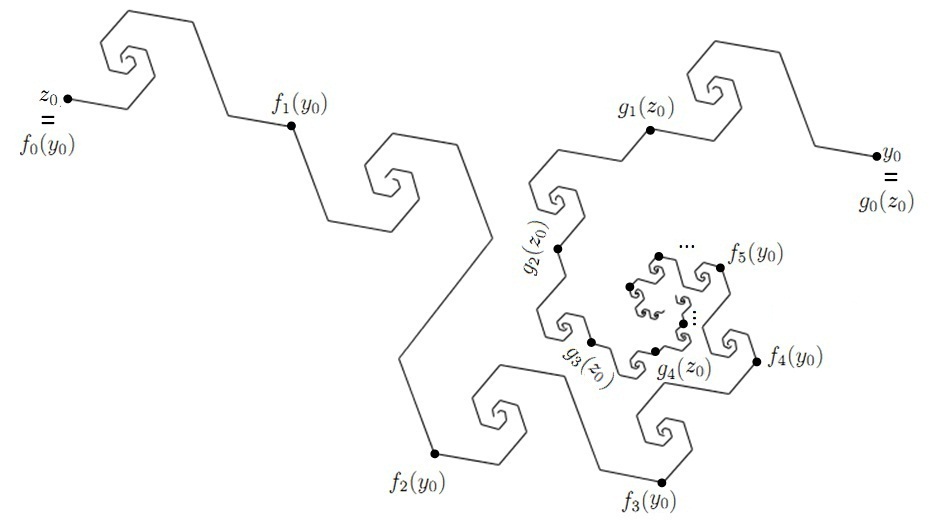}
\end{figure}
\centerline{Figure 9: Approximating $\mathcal{R}_{1-2\alpha}$ by $\varphi_{1}([0,1])$ and $\varphi_{2}([0,1])$.}

%\bigskip
%Let us zoom in the interval $[f_{0}(y_{0}),f_{1}(y_{0})]$ in $\varphi_{2}[0,1]$:

\begin{center}
    \includegraphics[scale=0.6]{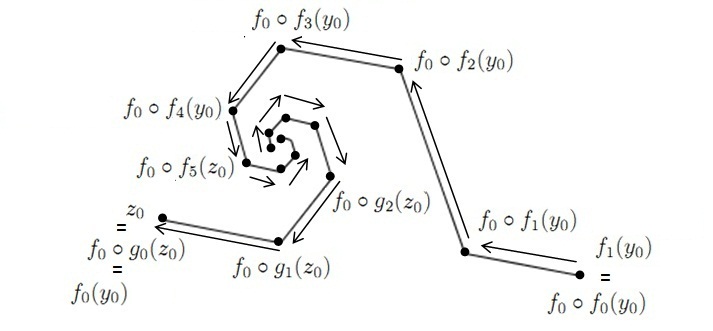}
  \end{center}
\centerline{Figure 10: Zoom of the interval $[f_{0}(y_{0}),f_{1}(y_{0})]$ in $\varphi_{2}([0,1])$. }

\bigskip

\begin{center}
    \includegraphics[width=3.8cm,angle=0]{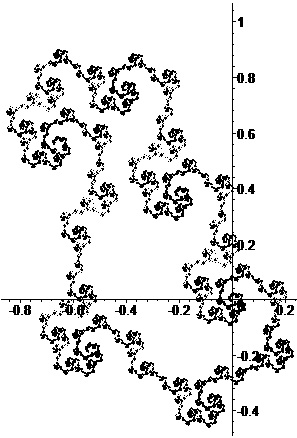}
  \end{center}
\centerline{Figure 11: Boundary of $\mathcal{R}_{3,-2}$.}

\bigskip

\noindent \textbf{Aknowledgments} I deeply thank Ali Messaoudi for all the help during the writing of the manuscript. \\

\bibliographystyle{plain}
\bibliography{Rauzy2016}

\textsc{Gustavo Antonio Pavani}

\textsc{Department of Mathematics - S\~{a}o Paulo State University - UNESP}

\textsc{Address: Rua Crist\'{o}v\~{a}o Colombo, 2265, Jardim Nazareth, 15054-000, S\~{a}o Jos\'{e} do Rio Preto, SP, Brazil}

\textsc{E-mail address:} \texttt{pavani@ibilce.unesp.br}

\end{document}